\documentclass[pdftex,11pt]{article}

\usepackage{amsmath}
\usepackage{amsfonts}
\usepackage{amssymb}

\usepackage{url}
\usepackage{graphicx}
\usepackage[usenames]{color}
\usepackage{subfigure}
\usepackage{latexsym}
\usepackage{hyperref}
\usepackage{pslatex}
\usepackage{verbatim}

\definecolor{DarkOliveGreen}{rgb}{0.33,0.42,0.18}
\definecolor{DodgerBlue4}{rgb}{0.06,0.31,0.55}

\hypersetup{
pdftitle={The Problem with the Linpack Benchmark 1.0 Matrix Generator},
pdfauthor={Julien Langou},
linkbordercolor={0 0 0},
citebordercolor={0 0 0},
urlbordercolor={0 0 0},
pdfborder={0 0 0},
colorlinks={true},
linkcolor=DodgerBlue4,
urlcolor=none,
citecolor=DodgerBlue4
}

\begin{document}

\noindent{\color{DodgerBlue4}\Large \textbf{
The Problem with the Linpack Benchmark 1.0 Matrix Generator}}\\

\noindent \textbf{\color{DodgerBlue4}
\lbrack v1\rbrack~Thu Jun 12, 2008.\\
\lbrack v2\rbrack~Thu Sep 18, 2008 (this version).
}\\

\noindent
\textbf{\color{DodgerBlue4}
Jack Dongarra}\\
Department of Electrical Engineering and Computer Science, University of Tennessee\\
Oak Ridge National Laboratory\\
University of Manchester\\
\textbf{\color{DodgerBlue4}
Julien Langou}\\
Department of Mathematical and Statistical Sciences, University of Colorado Denver\\

\begin{center}
\begin{minipage}{14cm}
\noindent\textbf{\color{DodgerBlue4}Abstract:} {\small

We characterize the matrix sizes for which the Linpack Benchmark 1.0 matrix
generator constructs a matrix with identical columns.

}
\end{minipage}
\end{center}

\section*{\color{DodgerBlue4}1~~~~Introduction}

Since 1993, twice a year, a list of the sites operating the 500 most powerful
computer systems is released by the TOP500 project~\cite{TOP500}. A single
number is used to rank computer systems based on the results obtained on the
\textit{High Performance Linpack Benchmark}~(HPL Benchmark).

The HPL Benchmark consists of solving a dense linear system in double
precision, 64--bit floating point arithmetic, using Gaussian elimination with
partial pivoting. The ground rules for running the benchmark state that the
supplied matrix generator, which uses a pseudo--random number generator, must
be used in running the HPL benchmark. The supplied matrix generator can be
found in \textit{High Performance Linpack 1.0}~(HPL--1.0)~\cite{HPL} which is an
implementation of the HPL Benchmark. In a HPL benchmark program, the
correctness of the computed solution is established and the performance is
reported in floating point operations per sec (flops/sec). It is this number
that is used to rank computer systems across the world in the TOP500 list.  For
more on the history and motivation for the HPL Benchmark, see~\cite{dolp:03}.

In May 2007, a large high performance computer manufacturer ran a
twenty-hour-long HPL Benchmark. The run fails with the output result:
\begin{verbatim}
   || A x - b ||_oo / ( eps * ||A||_1  * N ) = 9.224e+94 ...... FAILED
\end{verbatim}
It turned out that the manufacturer chose $n$ to be $n=2,220,032 = 2^{13}\cdot
271$.  This was a bad choice. In this case, the HPL Benchmark 1.0 matrix
generator produced a matrix $A$ with identical columns. Therefore the matrix
used in the test was singular and one of the checks of correctness determined
that there was a problem with the solution and the results should be considered
questionable. The reason for the suspicious results was neither a hardware
failure nor a software failure but a predictable numerical issue.

Nick Higham pointed out that this numerical issue had already been detected
in 1989 for the LINPACK-D benchmark implementation, a predecessor of HPL, and
had been reported to the community by David Hough~\cite{houg:89}. Another
report has been made to the HPL developers in 2004 by David Bauer with
$n=131,072$.  In this manuscript, we explain why and when the Linpack Benchmark 1.0
matrix generator generates matrices with identical columns.  We define
$\mathcal{S}$ as the set of all integers such that the Linpack Benchmark 1.0 matrix
generator produces a matrix with at least two identical columns.  We
characterize and give a simple algorithm to determine if a given $n$ is in
$\mathcal{S}$.\\

\noindent\textbf{ \color{DodgerBlue4} Definition~~1}
\textit{
We define $\mathcal{S}$ as the set of all integers such that the Linpack
Benchmark 1.0 matrix generator produces a matrix with at least two identical
columns. For $i > 2$, we define $\mathcal{S}_i$ as the set of all integers
such that the Linpack Benchmark 1.0 matrix generator produces a matrix with at
least one column repeated $i$ times.\\
}\\

In Table~\ref{tab:1}, for illustration, we give the $40$ smallest integers in
$\mathcal{S}$ along with the largest $i$ for which the associated matrix size is in
$\mathcal{S}_i$.

\begin{table}
\begin{center}
\begin{tabular}{rrrrr}
    65,536 (~2) &  98,304 (~2) & 131,072 (~8) & 147,456 (~2) & 163,840 (~3) \\
   180,224 (~2) & 196,608 (~6) & 212,992 (~2) & 229,376 (~4) & 245,760 (~2) \\
   262,144 (32) & 270,336 (~2) & 278,528 (~3) & 286,720 (~2) & 294,912 (~5) \\
   303,104 (~2) & 311,296 (~3) & 319,488 (~2) & 327,680 (10) & 335,872 (~2) \\
   344,064 (~3) & 352,256 (~2) & 360,448 (~6) & 368,640 (~2) & 376,832 (~3) \\
   385,024 (~2) & 393,216 (24) & 401,408 (~2) & 409,600 (~4) & 417,792 (~2) \\
   425,984 (~7) & 434,176 (~2) & 442,368 (~4) & 450,560 (~2) & 458,752 (14) \\
   466,944 (~2) & 475,136 (~4) & 483,328 (~2) & 491,520 (~8) & 499,712 (~2) \\
\end{tabular}
\end{center}
\color{DodgerBlue4}
\caption{\label{tab:1}
\color{black}
The $40$ matrix sizes smaller than $500,000$ for which the Linpack Benchmark 1.0
matrix generator will produce a matrix with identical columns. The number in
parenthesis indicates the maximum of the number of times each column is repeated.
For example, the entry ``491,520 (~8)'' indicates that, for the matrix size 491,520,
there exists one column that is repeated eight times while there exists no column that 
is repeated nine times.
}
\end{table}

Some remarks are in order.
\begin{enumerate}

\item[] \textbf{\color{DodgerBlue4}Remark~~1.1~~}
If $i > j > 2$ then $\mathcal{S}_i \subset \mathcal{S}_j \subset \mathcal{S}$.

\item[] \textbf{\color{DodgerBlue4}Remark~~1.2~~}
If $n$ is in $\mathcal{S}$, then the matrix generated by the Linpack
Benchmark 1.0 matrix generator has at least two identical columns, therefore this matrix is
necessarily singular.  If $n$ is not in $\mathcal{S}$, the coefficient matrix
has no identical columns; however we do not claim that the matrix is
nonsingular. Not being in $\mathcal{S}$ is not a sufficient condition for
being nonsingular.

\item[] \textbf{\color{DodgerBlue4}Remark~~1.3~~}
In practice, we would like the coefficient matrix to be well-conditioned
(since we want to numerically solve a linear system of equations associated
with them). This is a stronger condition than being nonsingular.  Edelman
in~\cite{edel:88} proves that for real $n$--by--$n$ matrices with elements from
a standard normal distribution, the expected value of the log of the 2-norm
condition number is asymptotic to $\log n$ as $n \to \infty$ (roughly $\log n +
1.537$).  The Linpack Benchmark 1.0 matrix generator uses a
uniform distribution on the interval \lbrack -0.5, 0.5\rbrack, for which the
expected value of the log of the 2-norm condition number is also asymptotic to
$\log n$ as $n \to \infty$ (roughly $4\log n + 1$), see Cuesta-Albertos and
Wschebor~\cite{caws:03}.  Random matrices are expected to be
well-conditioned; however, pseudo random number generator are only an attempt
to create randomness and we will see that, in some particular cases, the
generated matrices have repeated columns and are therefore singular (that is to
say infinitely ill--conditioned).

\item[] \textbf{\color{DodgerBlue4}Remark~~1.4~~}
HPL--1.0 checks whether a zero-pivot occurs during the factorization and
reports it to the user.  Due to rounding errors, even if the initial matrix has
two identical columns, exact-zero pivots hardly ever occur in practice.
Consequently, it is difficult for benchmarkers to distinguish between numerical
failures and hardware/software failures. This issue is further investigated in
\color{DodgerBlue4}\S5\color{black}.

\item[] \textbf{\color{DodgerBlue4}Remark~~1.5~~}
In Remark~\color{DodgerBlue4}1.3\color{black}, we stated that we would like the
coefficient matrix to be well-conditioned. Curiously enough, we will see in
\color{DodgerBlue4}\S5\color{black}  that the HPL benchmark can successfully return
when ran on a matrix with several identical columns. This is due to the
combined effect of finite precision arithmetic (that transforms a singular
matrix into an ill--conditioned matrix) and the use of a test for correctness
that is independent of the condition number of the coefficient matrix.

\end{enumerate}

\section*{\color{DodgerBlue4}2~~~~How the Linpack Benchmark matrix generator constructs a pseudo--random matrix}

The pseudo--random coefficient matrix $A$ from the HPL Benchmark 1.0 matrix
generator is generated by the HPL subroutine \texttt{HPL\_pdmatgen.c}. In this
subroutine, the pseudo--random number generator uses a linear congruential
algorithm~(see for example \cite[\S3.2]{knuth:97}) $$X(n+1) = (a * X(n) + c)
\textmd{ mod } m,$$ with $m=2^{31}$, $a=1103515245$, $c=1235$.  These choices
of $m$, $a$ and $c$ are fairly standard and we find them for example in the
standard POSIX.1-2001 or in the GNU libc library for the \texttt{rand()}
function.  The maximum period of a sequence generated by a linear congruential
algorithm is at most $m$, and in our case, with HPL--1.0's parameters $a$ and
$c$, we indeed obtain the maximal period $2^{31}$. (Proof: either by direct
check or using the Full-Period Theorem, see~\cite[\S3.2]{knuth:97}).
This provides us with a periodic sequence $s$ such that $s(i+2^{31}) = s(i),
\textmd{ for any } i\in \mathbb{N}$. HPL--1.0 fills its matrices with
pseudo--random numbers by columns using this sequence $s$ starting with $A(1,1)
= s(1)$, $A(2,1) = s(2)$, $A(3,1) = s(3)$, and so on.\\

\noindent\textbf{ \color{DodgerBlue4} Definition~~2}
\textit{
We define a Linpack Benchmark 1.0 matrix generator, a matrix generator such that
\begin{equation}
\label{hyp:1} 
A( i , j ) = s( ( j - 1 ) * n + i ),\quad 1\leq i,j\leq n.
\end{equation}
and $s$ is such that
\begin{equation}
\label{hyp:2} 
s(i+2^{31}) = s(i),\quad \textmd{ for any } i\in \mathbb{N}
\quad \textmd{ and } \quad
s(i) \ne s(j),\quad \textmd{ for any } 1\leq i,j \leq 2^{31}.
\end{equation}
}

Some remarks:
\begin{enumerate}

\item[] \textbf{\color{DodgerBlue4}Remark~~2.1~~}
The assumption $s(i) \ne s(j)$, for any $1\leq i,j \leq 2^{31}$ is true in the
case of the Linpack Benchmark 1.0 matrix generator. It can be relaxed to admit more
sequences $s$ for which some elements can be identical. However this assumption
makes the sufficiency proof of the theorem in \S{\color{DodgerBlue4}4} easier
and clearer.

\item[] \textbf{\color{DodgerBlue4}Remark~~2.2~~}
It is important to note that the matrix generated by the Linpack Benchmark 1.0
matrix generator solely depends on the dimension $n$.  The Linpack Benchmark 1.0
matrix generator requires benchmarkers to use the same matrix for any block
size, for any number of processors or for any grid size.

\item[] \textbf{\color{DodgerBlue4}Remark~~2.3~~}
Moreover, since the Linpack Benchmark 1.0 matrix generator possesses its own
implementation of the pseudo--random number generator, the computed pseudo--random
numbers in the sequence $s$ depend weakly on the computer systems.
Consequently the pivot pattern of the Gaussian elimination is preserved from
one computer system to the other, from one year to the other.

\item[] \textbf{\color{DodgerBlue4}Remark~~2.4~~}
Finally, the linear congruential algorithm for the sequence $s$ enables the
matrix generator for a scalable implementation of the construction of the
matrix: each process can generate their local part of the global matrix without
communicating or generating the global matrix. This property is not usual among
pseudo--random number generators.

\item[] \textbf{\color{DodgerBlue4}Remark~~2.5~~}
To give a sense of the magnitude of the size $n$ of matrices, the matrix size
for the \#1 entry in the TOP500 list of June 2008 was $2,236,927$ which is
between $2^{21}$ and $2^{22}$. The smallest matrix size in the TOP 500 list of
June 2008 was $273,919$ which is between $2^{18}$ and $2^{19}$.

\item[] \textbf{\color{DodgerBlue4}Remark~~2.6~~}
The pseudo--random number generator has been changed five times in the history 
of the Linpack Benchmark. We recall here some historical facts.
\begin{description}

\item[1980 -- \color{DodgerBlue4}LINPACKD--1.0 \color{black} --]
The initial LINPACKD benchmark uses a matrix generator based on the (Fortran) code below:\\
\includegraphics[width=0.50 \textwidth]{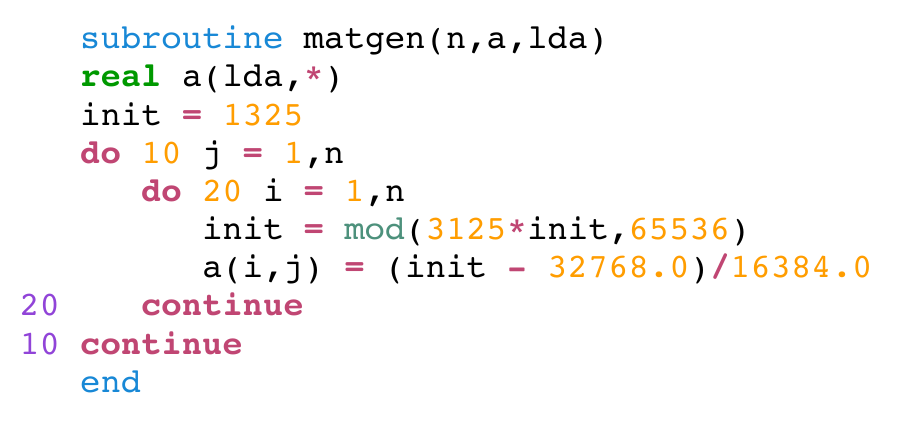}\\
The period of this pseudo--random number generator is:  $2^{14} = 16,384$.

\item[1989 -- numerical failure report --] David Hough~\cite{houg:89} observed
a numerical failure with the LINPACKD--1.0 benchmark for a matrix size $n=512$
and submitted his problem as an open question to the community through
NA-Digest.

\item[1989 -- \color{DodgerBlue4}LINPACKD--2.0 \color{black} --]
Two weeks after David Hough's post, Robert Schreiber~\cite{schr:89} posted in NA
Digest an explanation of the problem, he gave credit to Nick
Higham and himself for the explanation. The problem \#27.4
in Nick Higham's \textit{Accuracy and Stability of Numerical Algorithms}
book~\cite{higham:2002:ASN} is inspired from this story.
Higham and Schreiber also provide a patch to improve the pseudo--random number generator. 
Replacing line 6 of the previous code\\
\includegraphics[width=0.40 \textwidth]{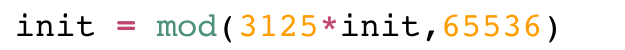}\\
by\\
\includegraphics[width=0.40 \textwidth]{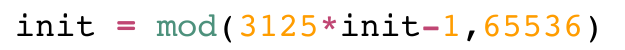}\\
increases the period from $2^{14} = 16,384$ to $2^{16} = 65,536$.
We call this version LINPACKD--2.0.

\item[1992 -- \color{DodgerBlue4}LINPACKD--3.0 \color{black} --] The
pseudo--random number generator of LINPACKD is updated for good in 1992 by
using the DLARUV LAPACK routine based on Fishman's multiplicative congruential
method with modulus $2^{48}$ and multiplier $33952834046453$
(see~\cite{fish:90}).

\item[2000 -- \color{DodgerBlue4}HPL--1.0 \color{black} --] First release of HPL (09/09/2000).
The pseudo--random number generator uses a linear congruential
algorithm~(see for example \cite[\S3.2]{knuth:97})
$$X(n+1) = (a * X(n) + c) \textmd{ mod } m,$$
with $m=2^{31}$, $a=1103515245$, $c=1235$.
The period of this pseudo--random number generator is $2^{31}$. 

\item[2004 -- numerical failure report --]
Gregory Bauer observed a numerical failure with HPL and $n=2^{17}=131,072$.
History repeats itself. The HPL developers recommended to HPL users willing to
test matrices of size larger than $2^{15}$ to not use power two.

\item[2007 -- numerical failure report --] A large manufacturer observed a numerical
failure with HPL and $n=2,220,032$. History repeats itself again. Note that
$2,200,032 = 2^{13}\cdot271$, and is not a power of two.

\item[2008 -- \color{DodgerBlue4}HPL--2.0 \color{black} --]
This present manuscript explains the problem in the Linpack Benchmark 1.0
matrix generator. As of September 10th 2008,  Piotr Luszczek has incorporated a
new pseudo--random number generator in HPL--2.0.  This pseudo--random number generator
uses a linear congruential algorithm with $a = 6364136223846793005$, $c = 11$
and $m = 2^{64}$.
The period of this pseudo--random number generator is $2^{64}$. 

\end{description}

\end{enumerate}

\section*{\color{DodgerBlue4}3~~~~Understanding~$\mathcal{S}$}

Consider a large dense matrix of order $3\cdot10^6$ generated by the process
described in Definition~{\color{DodgerBlue4}2}. The number of entries in this
matrix is $9\cdot10^{12}$ which is above the pseudo--random number generator period
($2^{31}\approx2.14\cdot10^{9}$). However, despite this fact, it is fairly
likely for the constructed matrix to have distinct columns and even to be
well--conditioned.

On the other hand, we can easily generate a ``small'' matrix with identical
columns. Take n=$2^{16}$, we have for any $i = 1,\ldots, n$:
   $$A(i,2^{15}+1) = s( i + n*(j-1) ) = s( i + 2^{15}*n ) = s(i + 2^{15} * 2^{16} ) = s(i+2^{31}) = s(i) = A(i,1),$$
therefore the column $1$ and the column $2^{15}+1$ are exactly the same. The
column $2$ and the column $2^{15}+2$ are exactly the same, etc.  We can
actually prove that $2^{16} = 65,536$ is the smallest matrix order for which a
multiple of a column can happen.

Another example of $n\in\mathcal{S}$ is $n=2^{31}=2,147,483,648$ for which all
columns of the generated matrix are the same.  Our goal in this section is to
build more $n$ in $\mathcal{S}$ to have a better knowledge of this set.\\

\textbf{\color{DodgerBlue4} If $n$ is a multiple of $2^0=1$ and $n>2^{31}$ then $n\in\mathcal{S}$.}
(Note that the statement ``any $n$ is multiple of $2^0=1$ and $n>2^{31}$'' means $n>2^{31}$.)
The reasoning is as follows. There are $2^{31}$ indexes from $1$ to $2^{31}$.
Since there are at least $2^{31}+1$ elements in the first row of $A$ (assumption $n>2^{31}$), then,
necessarily, at least one index (say $k$) is repeated twice in the first row of $A$.
This is the pigeonhole principle.
Therefore we have proved the existence of two columns $i$ and $j$ such that they both
start with the $k$--th term of the sequence. If two columns start with the index of the sequence,
they are the same (since we take the element of the column sequentially in the sequence).
The three smallest numbers of this type are
\begin{eqnarray}
\nonumber n = 2^{0} * (2^{31} + 1) = 2,147,483,649 \in \mathcal{S}\\
\nonumber n = 2^{0} * (2^{31} + 2) = 2,147,483,650 \in \mathcal{S}\\
\nonumber n = 2^{0} * (2^{31} + 3) = 2,147,483,651 \in \mathcal{S}
\end{eqnarray}

\textbf{\color{DodgerBlue4} If $n$ is a multiple of $2^1=2$ and $n>2^{30}$ then $n\in\mathcal{S}$.}
If $n$ is even ($n = 2q$), then the first row of $A$ accesses the numbers of
the sequence $s$ using only odd indexes. There are $2^{30}$ odd indexes
between $1$ and $2^{31}$.
Since there are at least $2^{30}+1$ elements in the first row of $A$ (assumption $n>2^{30}$), then,
necessarily, at least one index is repeated twice in the first row of $A$.
This is the pigeonhole principle.
The three smallest numbers of this type are:
\begin{eqnarray}
\nonumber n = 2^{1} * (2^{29} + 1) = 1,073,741,826 \in \mathcal{S}\\
\nonumber n = 2^{1} * (2^{29} + 2) = 1,073,741,828 \in \mathcal{S}\\
\nonumber n = 2^{1} * (2^{29} + 3) = 1,073,741,830 \in \mathcal{S}.
\end{eqnarray}

\textbf{\color{DodgerBlue4} If $n$ is a multiple of $2^2=4$ and $n>2^{29}$ then $n\in\mathcal{S}$.}
If $n$ is a multiple of $4$ ($n = 4q$), then the first row of $A$ accesses the numbers of
the sequence $s$ using only $(4q+1)$--indexes. There are $2^{29}$ $(4q+1)$--indexes
between $1$ and $2^{31}$.
Since there are at least $2^{29}+1$ elements in the first row of $A$ (assumption $n>2^{29}$), then,
necessarily, at least one index is repeated twice in the first row of $A$.
This is the pigeonhole principle.
The first three numbers of this type are:
\begin{eqnarray}
\nonumber n = 2^{2} * (2^{27} + 1) = 536,870,916 \in \mathcal{S}\\
\nonumber n = 2^{2} * (2^{27} + 2) = 536,870,920 \in \mathcal{S}\\
\nonumber n = 2^{2} * (2^{27} + 3) = 536,870,924 \in \mathcal{S}.
\end{eqnarray}
$$\vdots$$

\textbf{\color{DodgerBlue4} If $n$ is a multiple of $2^{13}$ and $n>2^{18}$ then $n\in\mathcal{S}$.}
This gives for example:
\begin{eqnarray}
\nonumber n_{12} = 2^{13} * (2^{5} + 1)  = 2^{13} * 33 & = & 270,336 \in \mathcal{S}\\
\nonumber n_{13} = 2^{13} * (2^{5} + 2)  = 2^{13} * 34 & = & 278,528 \in \mathcal{S}\\
\nonumber n_{15} = 2^{13} * (2^{5} + 3)  = 2^{13} * 35 & = & 294,912 \in \mathcal{S}.
\end{eqnarray}
These three numbers correspond to entries $(3,2)$, $(3,3)$ and $(3,5)$ in Table~\ref{tab:1}.\\

\textbf{\color{DodgerBlue4} If $n$ is a multiple of $2^{14}$ and $n>2^{17}$ then $n\in\mathcal{S}$.}
This gives for example:
\begin{eqnarray}
\nonumber n_4 = 2^{14} * (2^{3} + 1) = 2^{14} *  9 & = & 147,456 \in \mathcal{S}\\
\nonumber n_5 = 2^{14} * (2^{3} + 2) = 2^{14} * 10 & = & 163,840 \in \mathcal{S}\\
\nonumber n_6 = 2^{14} * (2^{3} + 3) = 2^{14} * 11 & = & 180,224 \in \mathcal{S}.
\end{eqnarray}
These three numbers correspond to entries $(1,4)$, $(1,5)$ and $(2,1)$ in Table~\ref{tab:1}.\\

\textbf{\color{DodgerBlue4} If $n$ is a multiple of $2^{15}$ and $n>2^{16}$ then $n\in\mathcal{S}$.}
This gives for example:
\begin{eqnarray}
\nonumber n_2 = 2^{15} * (2^{1} + 1) = 2^{15} *  3 & = &  98,304 \in \mathcal{S}\\
\nonumber n_3 = 2^{15} * (2^{1} + 2) = 2^{15} *  4 & = & 131,072 \in \mathcal{S}\\
\nonumber n_5 = 2^{15} * (2^{1} + 3) = 2^{15} *  5 & = & 163,840 \in \mathcal{S}.
\end{eqnarray}
These three numbers correspond to entries $(1,2)$, $(1,3)$ and $(1,5)$ in Table~\ref{tab:1}.\\

\textbf{\color{DodgerBlue4} If $n$ is a multiple of $2^{16}$ and $n>2^{15}$ then $n\in\mathcal{S}$.}
\begin{eqnarray}
\nonumber n_1 = 2^{16} * (2^{0} + 1) = 2^{16} *  1 & = &  65,536 \in \mathcal{S}\\
\nonumber n_3 = 2^{16} * (2^{0} + 2) = 2^{16} *  2 & = & 131,072 \in \mathcal{S}\\
\nonumber n_7 = 2^{16} * (2^{0} + 3) = 2^{16} *  3 & = & 196,608 \in \mathcal{S}.
\end{eqnarray}
These three numbers correspond to entries $(1,1)$, $(1,3)$ and $(2,2)$ in Table~\ref{tab:1}.

From this section, we understand that any $n$ multiple of $2^k$ and larger than
$2^{31-k}$ is in $\mathcal{S}$.  In the next paragraph, we prove that this is
indeed the only integers in $\mathcal{S}$ which provides us with a complete
characterization of $\mathcal{S}$.

\section*{\color{DodgerBlue4}4~~~~Characterization of $\mathcal{S}$}

\textbf{Theorem:}
$n \in \mathcal{S}$
if and only if
the matrix of size $n$ generated by the Linpack Benchmark 1.0 matrix generator has at least two identical columns
if and only if
$$ n > 2^{31 - k } \quad \textmd{ where } n = 2^k\cdot q \textmd{ with } q \textmd{ odd}.  $$

\textit{Proof:}

\begin{itemize}
\item[\begin{tabular}{|c|}\hline$\Leftarrow$\\\hline\end{tabular}]
Let us assume that $n$ is a multiple of $2^k$, that is to say
$$n = 2^k\cdot q, \quad 1 \leq q $$
and let us assume that
$$ n > 2^{31 - k }.$$
In this case, the first row of $A$ accesses the numbers of
the sequence $s$ using only $( 2^k \cdot q + 1 )$--indexes.
There are $2^{31-k}$ $( 2^k \cdot q + 1 )$--indexes
between $1$ and $2^{31}$.
Since there are at least $2^{31-k}+1$ elements in the first row of $A$ (assumption $n>2^{31-k}$), then,
necessarily, at least one index is repeated twice in the first row of $A$.
This is the pigeonhole principle.
If two columns start with the same index in the sequence,
they are the same (since we take the element of the column sequentially in the sequence).

\item[\begin{tabular}{|c|}\hline$\Rightarrow$\\\hline\end{tabular}]
Assume that there are two identical columns $i$ and $j$ in the
matrix generated by the Linpack Benchmark 1.0 matrix generator ($i \ne j$). Without loss of generality, assume $i>j$.
The fact that column $i$ is the same as column $j$ means that these columns have identical entries, in particular,
they share the same first entry. We have
$$ A(1,i) = A(1,j). $$
From this, Equation~(\ref{hyp:1}) implies
$$ s\left(1+(i-1)n\right) = s\left(1+(j-1)n\right). $$
Equation~(\ref{hyp:2}) states that all elements in a period of length $2^{31}$ are different,
therefore, since $i\ne j$, we necessarily have
$$ 1+(i-1)n = 1+(j-1)n + 2^{31} \cdot p, \quad 1 \leq p. $$
This implies
$$ (i-j) n = 2^{31} \cdot p, \quad 1 \leq p. $$
We now use the fact that $ n= 2^k \cdot q$ with $q$ odd and get
$$ (i-j) \cdot 2^k \cdot q = 2^{31} \cdot p, \quad 1 \leq p, \quad q \textmd{ is odd}. $$
Since $q$ is odd, this last equality implies that $2^{31}$ is a divisor of $(i-j)\cdot2^k$. This writes
$$ (i-j)\cdot2^k =  2^{31} \cdot r, \quad 1 \leq r.$$
From which, we deduce that
$$ (i-j)\cdot2^k \geq  2^{31} .$$
A upper bound for $i$ is $n$, a lower bound for $j$ is $1$; therefore,
$$ (n-1)\cdot2^k \geq  2^{31} .$$

We conclude that, if a matrix of size $n$ generated by the Linpack Benchmark 1.0 matrix generator has at least two identical columns, this implies
$$ n > 2^{31 - k } \quad \textmd{ where } n = 2^k\cdot q \textmd{ with } q \textmd{ odd}.  $$

\hfill$\square$
\end{itemize}

\section*{\color{DodgerBlue4}5~~~~Solving (exactly) singular system in finite precision arithmetic
with a small backward error}

From our analysis, the first matrix size $n$ for which the Linpack Benchmark 1.0 matrix generator 
will generate a matrix with two identical columns is $n = 65,536$
(see Table~\color{DodgerBlue4}1\color{black}). However, HPL--1.0 passes
all the test for correctness on this matrix size. The same for $n = 98,304$ which is our second
matrix size in the list
(see Table~\color{DodgerBlue4}1\color{black}). If we look more carefully at the output file for $n=2,220,032$, we see
that only one out of the three test for correctness is triggered:
\begin{verbatim}
   ||Ax-b||_oo / ( eps * ||A||_1  * N        ) =        9.224e+94 ...... FAILED
   ||Ax-b||_oo / ( eps * ||A||_1  * ||x||_1  ) =        0.0044958 ...... PASSED
   ||Ax-b||_oo / ( eps * ||A||_oo * ||x||_oo ) =        0.0000002 ...... PASSED
\end{verbatim}
Despite the fact that the matrix has identical columns, we observe that
HPL--1.0 is enable to pass sometimes all the tests, sometimes two tests out of
three, sometimes none of the three tests. This section will answer how this
behavior is possible. First of all, we need to explain how the Linpack Benchmark assesses
the correctness of an answer.

\subsection*{\color{DodgerBlue4}5.1~~~~How the Linpack Benchmark program checks a solution}

To verify the result after the LU factorization, the benchmark regenerates the input matrix and the
right-hand side, then an accuracy check on the residual $Ax-b$ is performed.

The LINPACKD benchmark checks the accuracy of the solution by returning
\begin{equation}
\nonumber
 \frac{\|Ax-b\|_{\infty}}{n\epsilon\|A\|_M\|x\|_\infty}
\end{equation}
where $\| A \|_M = \max_{i,j} |a_{ij} |$.  and $\epsilon$ is the relative
machine precision.

For HPL--1.0, the three following scaled residuals are computed
\begin{equation}
\nonumber
 r_n = \frac{\|Ax-b\|_{\infty}}{n\epsilon\|A\|_1},
\end{equation}
\begin{equation}
\nonumber
 r_1 = \frac{\|Ax-b\|_{\infty}}{\epsilon\|A\|_1 \|x\|_1},
\end{equation}
\begin{equation}
\nonumber
 r_{\infty} = \frac{\|Ax-b\|_{\infty}}{n\epsilon\|A\|_{\infty}\|x\|_{\infty}}.
\end{equation}
A solution is considered numerically correct when all of these quantities are
less than a threshold value of 16.  The last quantity ( $r_\infty$ )
corresponds to the normwise backward error in the infinite norm allowing
perturbations on $A$ only~\cite{higham:2002:ASN}. The last two quantities (~$r_\infty$, $r_1$~) are
independent of the condition number of the coefficient matrix $A$ and should
always be less than a threshold value of the order of $1$ (no matter how
ill--conditioned $A$ is).

As of HPL--2.0, the check for correctness is
\begin{equation}
r_4 = \frac{\|Ax-b\|_{\infty}}{n\epsilon\left(\|A\|_{\infty}\|x\|_{\infty}+\|b\|_{\infty}\right)}.
\end{equation}
This corresponds to the normwise backward error in the infinite norm allowing
perturbations on $A$ and $b$ only~\cite{higham:2002:ASN}.  A solution is considered numerically correct
when this quantity is less than a threshold value of 16. Although the error
analysis of Gaussian elimination with partial pivoting can be done in such a
way that $b$ is not perturbed (in other words $r_{\infty}$ is the criterion you
want to use for Gaussian elimination with partial pivoting), HPL--2.0 switches
to $r_4$, the usual backward error as found in textbooks.

This discussion on the check for correctness explains why HPL--1.0 is able to
pass the test for correctness even though the input matrix is exactly singular.

\subsection*{\color{DodgerBlue4}5.2~~~~Repeating identical blocks to the underflow}

In~\cite{schr:89}, Schreiber and Higham explain what happens when a block is
repeated $k$ times in the initial coefficient matrix $A$. At each repeat, the
magnitude of the pivot (diagonal entries of the $U$ matrix) are divided by
$\varepsilon$.  This is illustrated in
Figure~\color{DodgerBlue4}1\color{black}.  This process continues until
underflow. Denormalized might help but the process is still the same and ultimately
a zero pivot is reached, and the algorithm is stopped. In single precision arithmetic with $\varepsilon_s =
2^{-24}$ and underflow $2^{-126}$, five identical blocks will lead to
underflow. In double precision arithmetic with $\varepsilon = 2^{-16}$ and
underflow $2^{-1022}$, one will need 64 identical blocks. 

\begin{figure}[!h]
\color{DodgerBlue4}
\begin{center}
\includegraphics[width=0.50 \textwidth]{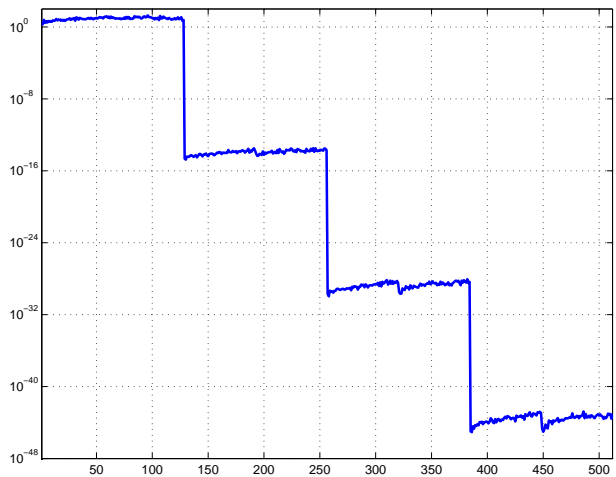}
\end{center}
\caption{\label{fig:}
\color{black}
Magnitude of the pivot (diagonal entries along the matrix $U$) for $n=512=2^9$ and 
the LINPACK--2.0 matrix generator.
The period of the LINPACK--2.0 matrix generator is $n=65536=2^{16}$ so that, for a matrix of size $n=512$,
columns repeat every 128 column. We observe that pivots are multiplied by $\varepsilon\approx 2.2\cdot10^{-16}$
at every repetition.
}
\end{figure}

\subsection*{\color{DodgerBlue4}5.3~~~~Anomalies in Matrix Sizes Reported in the June 2008 Top500 List}

Readers of this manuscript may be surprised to find three entries in the TOP
500 data from June 2008 with matrix sizes that lead to matrices with identical
columns if the HPL test matrix generator is used.  These three entries are
given in~Table~\ref{tab:2}.  For example, the run for the Earth Simulator from
2002 was done with $n = 1,075,200$ which corresponds to $2^{11} \cdot 525$,
therefore, the column $j = 2^{20} = 1,048,576$ would have been a repeat of the
first under our assumptions. The benchmark run on the Earth Simulator in 2002
was done with an older version of the test harness. This test harness predates
the HPL test harness and uses another matrix generator than the one provided by
HPL.  Today we require the HPL test harness to be used in the benchmark run.

\begin{table}
\begin{center}
\begin{tabular}{clllr}
Rank & Site & Manufacturer
& Year
& NMax \\
16 & Information Technology Center, The University of Tokyo & Hitachi
& 2008
& 1,433,600 (6) \\
49 & The Earth Simulator Center & NEC
& 2002
& 1,075,200 (2) \\
88 & Cardiff University - ARCCA & Bull SA
& 2008
& 634,880 (2) \\
\end{tabular}
\end{center}
\color{DodgerBlue4}
\caption{\label{tab:2}
\color{black}
The three entries in the TOP500 June 2008 list with suspicious $n$.
}
\end{table}

\section*{\color{DodgerBlue4}6~~~~How to fix the problem}

Between $1$ and $1\cdot10^6$, there are $49$ matrix sizes in $\mathcal{S}$ (see
Table~\ref{tab:1}).  Between $1$ and $3\cdot10^6$, there are $1,546$ matrix
sizes in $\mathcal{S}$ (see Appendix~{\color{DodgerBlue4}B}). Therefore, for this order of
matrix size, there is a good chance to choose a matrix size that is
not in $\mathcal{S}$. Unfortunately benchmarkers tend to pick multiples of high power of 2 for their matrix sizes
which increases the likelihood of picking an $n\in \mathcal{S}$.

\begin{enumerate}

\item The obvious recommendation is to choose any $n$ as long as it is odd.
In the odd case if $n < 2^{31} \approx 4\cdot 10^9 $, then $n \notin \mathcal{S}$.

\item A check can be added at the beginning of the execution of the Linpack Benchmark matrix generator.
The C-code looks as follows:

\includegraphics[width= 0.5\textwidth]{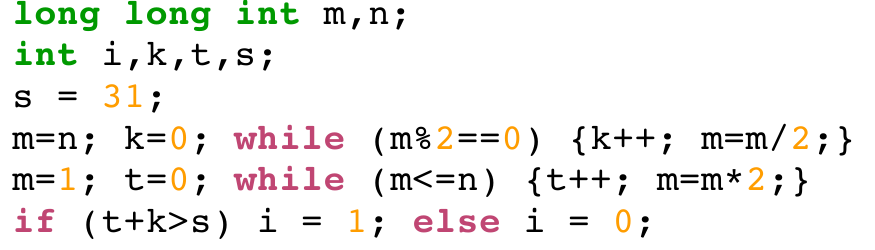}

$n$ is the matrix size, $2^s$ is period of the pseudo--random number generator ($s=31$
in our case) and $i$ is the output flag.  If $i = 1$, then $n\in\mathcal{S}$.
 If $i = 0$, then $n\notin\mathcal{S}$.
(The check could also consist of looking over the data given in Appendix~{\color{DodgerBlue4}B}).

\item If $n\in\mathcal{S}$, one can simply pad the matrix with an extra line. This can be easily done in
the HPL code \texttt{HPL\_pdmatgen.c} by changing the variable \texttt{jump3} from \texttt{M} to \texttt{M+1}
whenever $n\in\mathcal{S}$.

\item Another possibility is to increase the period of the pseudo--random number
generator used.  For example, if the pseudo--random number generator had a period of
$2^{64}$ and if $n \leq 2^{32}$,  then, assuming $( i\ne j \Rightarrow s(i)\ne s(j) )$,
entries would never repeat.

\item A check for correctness robust to ill--conditioned matrix could be used as discussed in \color{DodgerBlue4}\S5\color{black}.

\end{enumerate}

The Problem with the Linpack Benchmark 1.0 matrix generator is now corrected in
the Linpack Benchmark 2.0 Matrix Generator. The fix includes both proposition
4 (extend the period of the pseudo random generator) and proposition 5 (have 
a test for correctness robust to ill--conditioned matrices).

\section*{\color{DodgerBlue4}Acknowledgments}
The authors would like to thank Piotr Luszczek and Antoine Petitet for their
valuable comments on HPL, Nick Higham for making us aware of
David Hough \textit{Random Story}~\cite{houg:89} and his comments on the
backward error analysis of Gaussian elimination with partial pivoting, and finally Asim Yarkhan for one pertinent observation.

\bibliographystyle{plain}
\bibliography{hpl}

\newpage
\appendix

\section*{\color{DodgerBlue4}A~~~~The $1,564$ matrix sizes of $n$
from $1$ to $3,000,000$ for which the Linpack Benchmark matrix generator will
construct a matrix with identical columns}
\small
\begin{tabular}{|rr|rr|rr|rr|rr|}
    65536&  (  2) &    98304&  (  2) &   131072&  (  8) &   147456&  (  2) &   163840&  (  3) \\
   180224&  (  2) &   196608&  (  6) &   212992&  (  2) &   229376&  (  4) &   245760&  (  2) \\
   262144&  ( 32) &   270336&  (  2) &   278528&  (  3) &   286720&  (  2) &   294912&  (  5) \\
   303104&  (  2) &   311296&  (  3) &   319488&  (  2) &   327680&  ( 10) &   335872&  (  2) \\
   344064&  (  3) &   352256&  (  2) &   360448&  (  6) &   368640&  (  2) &   376832&  (  3) \\
   385024&  (  2) &   393216&  ( 24) &   401408&  (  2) &   409600&  (  4) &   417792&  (  2) \\
   425984&  (  7) &   434176&  (  2) &   442368&  (  4) &   450560&  (  2) &   458752&  ( 14) \\
   466944&  (  2) &   475136&  (  4) &   483328&  (  2) &   491520&  (  8) &   499712&  (  2) \\
   507904&  (  4) &   516096&  (  2) &   524288&  (128) &   528384&  (  2) &   532480&  (  3) \\
   536576&  (  2) &   540672&  (  5) &   544768&  (  2) &   548864&  (  3) &   552960&  (  2) \\
   557056&  (  9) &   561152&  (  2) &   565248&  (  3) &   569344&  (  2) &   573440&  (  5) \\
   577536&  (  2) &   581632&  (  3) &   585728&  (  2) &   589824&  ( 18) &   593920&  (  2) \\
   598016&  (  3) &   602112&  (  2) &   606208&  (  5) &   610304&  (  2) &   614400&  (  3) \\
   618496&  (  2) &   622592&  ( 10) &   626688&  (  2) &   630784&  (  3) &   634880&  (  2) \\
   638976&  (  5) &   643072&  (  2) &   647168&  (  3) &   651264&  (  2) &   655360&  ( 40) \\
   659456&  (  2) &   663552&  (  3) &   667648&  (  2) &   671744&  (  6) &   675840&  (  2) \\
   679936&  (  3) &   684032&  (  2) &   688128&  ( 11) &   692224&  (  2) &   696320&  (  3) \\
   700416&  (  2) &   704512&  (  6) &   708608&  (  2) &   712704&  (  3) &   716800&  (  2) \\
   720896&  ( 22) &   724992&  (  2) &   729088&  (  3) &   733184&  (  2) &   737280&  (  6) \\
   741376&  (  2) &   745472&  (  3) &   749568&  (  2) &   753664&  ( 12) &   757760&  (  2) \\
   761856&  (  3) &   765952&  (  2) &   770048&  (  6) &   774144&  (  2) &   778240&  (  3) \\
   782336&  (  2) &   786432&  ( 96) &   790528&  (  2) &   794624&  (  4) &   798720&  (  2) \\
   802816&  (  7) &   806912&  (  2) &   811008&  (  4) &   815104&  (  2) &   819200&  ( 13) \\
   823296&  (  2) &   827392&  (  4) &   831488&  (  2) &   835584&  (  7) &   839680&  (  2) \\
   843776&  (  4) &   847872&  (  2) &   851968&  ( 26) &   856064&  (  2) &   860160&  (  4) \\
   864256&  (  2) &   868352&  (  7) &   872448&  (  2) &   876544&  (  4) &   880640&  (  2) \\
   884736&  ( 14) &   888832&  (  2) &   892928&  (  4) &   897024&  (  2) &   901120&  (  7) \\
   905216&  (  2) &   909312&  (  4) &   913408&  (  2) &   917504&  ( 56) &   921600&  (  2) \\
   925696&  (  4) &   929792&  (  2) &   933888&  (  8) &   937984&  (  2) &   942080&  (  4) \\
   946176&  (  2) &   950272&  ( 15) &   954368&  (  2) &   958464&  (  4) &   962560&  (  2) \\
   966656&  (  8) &   970752&  (  2) &   974848&  (  4) &   978944&  (  2) &   983040&  ( 30) \\
   987136&  (  2) &   991232&  (  4) &   995328&  (  2) &   999424&  (  8) &  1003520&  (  2) \\
  1007616&  (  4) &  1011712&  (  2) &  1015808&  ( 16) &  1019904&  (  2) &  1024000&  (  4) \\
  1028096&  (  2) &  1032192&  (  8) &  1036288&  (  2) &  1040384&  (  4) &  1044480&  (  2) \\
  1048576&  (512) &  1050624&  (  2) &  1052672&  (  3) &  1054720&  (  2) &  1056768&  (  5) \\
  1058816&  (  2) &  1060864&  (  3) &  1062912&  (  2) &  1064960&  (  9) &  1067008&  (  2) \\
  1069056&  (  3) &  1071104&  (  2) &  1073152&  (  5) &  1075200&  (  2) &  1077248&  (  3) \\
  1079296&  (  2) &  1081344&  ( 17) &  1083392&  (  2) &  1085440&  (  3) &  1087488&  (  2) \\
  1089536&  (  5) &  1091584&  (  2) &  1093632&  (  3) &  1095680&  (  2) &  1097728&  (  9) \\
  1099776&  (  2) &  1101824&  (  3) &  1103872&  (  2) &  1105920&  (  5) &  1107968&  (  2) \\
  1110016&  (  3) &  1112064&  (  2) &  1114112&  ( 34) &  1116160&  (  2) &  1118208&  (  3) \\
  1120256&  (  2) &  1122304&  (  5) &  1124352&  (  2) &  1126400&  (  3) &  1128448&  (  2) \\
  1130496&  (  9) &  1132544&  (  2) &  1134592&  (  3) &  1136640&  (  2) &  1138688&  (  5) \\
  1140736&  (  2) &  1142784&  (  3) &  1144832&  (  2) &  1146880&  ( 18) &  1148928&  (  2) \\
  1150976&  (  3) &  1153024&  (  2) &  1155072&  (  5) &  1157120&  (  2) &  1159168&  (  3) \\
  1161216&  (  2) &  1163264&  (  9) &  1165312&  (  2) &  1167360&  (  3) &  1169408&  (  2) \\
  1171456&  (  5) &  1173504&  (  2) &  1175552&  (  3) &  1177600&  (  2) &  1179648&  ( 72) \\
  1181696&  (  2) &  1183744&  (  3) &  1185792&  (  2) &  1187840&  (  5) &  1189888&  (  2) \\
  1191936&  (  3) &  1193984&  (  2) &  1196032&  ( 10) &  1198080&  (  2) &  1200128&  (  3) \\
  1202176&  (  2) &  1204224&  (  5) &  1206272&  (  2) &  1208320&  (  3) &  1210368&  (  2) \\
  1212416&  ( 19) &  1214464&  (  2) &  1216512&  (  3) &  1218560&  (  2) &  1220608&  (  5) \\
\end{tabular}
\newpage
\small
\begin{tabular}{|rr|rr|rr|rr|rr|}
  1222656&  (  2) &  1224704&  (  3) &  1226752&  (  2) &  1228800&  ( 10) &  1230848&  (  2) \\
  1232896&  (  3) &  1234944&  (  2) &  1236992&  (  5) &  1239040&  (  2) &  1241088&  (  3) \\
  1243136&  (  2) &  1245184&  ( 38) &  1247232&  (  2) &  1249280&  (  3) &  1251328&  (  2) \\
  1253376&  (  5) &  1255424&  (  2) &  1257472&  (  3) &  1259520&  (  2) &  1261568&  ( 10) \\
  1263616&  (  2) &  1265664&  (  3) &  1267712&  (  2) &  1269760&  (  5) &  1271808&  (  2) \\
  1273856&  (  3) &  1275904&  (  2) &  1277952&  ( 20) &  1280000&  (  2) &  1282048&  (  3) \\
  1284096&  (  2) &  1286144&  (  5) &  1288192&  (  2) &  1290240&  (  3) &  1292288&  (  2) \\
  1294336&  ( 10) &  1296384&  (  2) &  1298432&  (  3) &  1300480&  (  2) &  1302528&  (  5) \\
  1304576&  (  2) &  1306624&  (  3) &  1308672&  (  2) &  1310720&  (160) &  1312768&  (  2) \\
  1314816&  (  3) &  1316864&  (  2) &  1318912&  (  6) &  1320960&  (  2) &  1323008&  (  3) \\
  1325056&  (  2) &  1327104&  ( 11) &  1329152&  (  2) &  1331200&  (  3) &  1333248&  (  2) \\
  1335296&  (  6) &  1337344&  (  2) &  1339392&  (  3) &  1341440&  (  2) &  1343488&  ( 21) \\
  1345536&  (  2) &  1347584&  (  3) &  1349632&  (  2) &  1351680&  (  6) &  1353728&  (  2) \\
  1355776&  (  3) &  1357824&  (  2) &  1359872&  ( 11) &  1361920&  (  2) &  1363968&  (  3) \\
  1366016&  (  2) &  1368064&  (  6) &  1370112&  (  2) &  1372160&  (  3) &  1374208&  (  2) \\
  1376256&  ( 42) &  1378304&  (  2) &  1380352&  (  3) &  1382400&  (  2) &  1384448&  (  6) \\
  1386496&  (  2) &  1388544&  (  3) &  1390592&  (  2) &  1392640&  ( 11) &  1394688&  (  2) \\
  1396736&  (  3) &  1398784&  (  2) &  1400832&  (  6) &  1402880&  (  2) &  1404928&  (  3) \\
  1406976&  (  2) &  1409024&  ( 22) &  1411072&  (  2) &  1413120&  (  3) &  1415168&  (  2) \\
  1417216&  (  6) &  1419264&  (  2) &  1421312&  (  3) &  1423360&  (  2) &  1425408&  ( 11) \\
  1427456&  (  2) &  1429504&  (  3) &  1431552&  (  2) &  1433600&  (  6) &  1435648&  (  2) \\
  1437696&  (  3) &  1439744&  (  2) &  1441792&  ( 88) &  1443840&  (  2) &  1445888&  (  3) \\
  1447936&  (  2) &  1449984&  (  6) &  1452032&  (  2) &  1454080&  (  3) &  1456128&  (  2) \\
  1458176&  ( 12) &  1460224&  (  2) &  1462272&  (  3) &  1464320&  (  2) &  1466368&  (  6) \\
  1468416&  (  2) &  1470464&  (  3) &  1472512&  (  2) &  1474560&  ( 23) &  1476608&  (  2) \\
  1478656&  (  3) &  1480704&  (  2) &  1482752&  (  6) &  1484800&  (  2) &  1486848&  (  3) \\
  1488896&  (  2) &  1490944&  ( 12) &  1492992&  (  2) &  1495040&  (  3) &  1497088&  (  2) \\
  1499136&  (  6) &  1501184&  (  2) &  1503232&  (  3) &  1505280&  (  2) &  1507328&  ( 46) \\
  1509376&  (  2) &  1511424&  (  3) &  1513472&  (  2) &  1515520&  (  6) &  1517568&  (  2) \\
  1519616&  (  3) &  1521664&  (  2) &  1523712&  ( 12) &  1525760&  (  2) &  1527808&  (  3) \\
  1529856&  (  2) &  1531904&  (  6) &  1533952&  (  2) &  1536000&  (  3) &  1538048&  (  2) \\
  1540096&  ( 24) &  1542144&  (  2) &  1544192&  (  3) &  1546240&  (  2) &  1548288&  (  6) \\
  1550336&  (  2) &  1552384&  (  3) &  1554432&  (  2) &  1556480&  ( 12) &  1558528&  (  2) \\
  1560576&  (  3) &  1562624&  (  2) &  1564672&  (  6) &  1566720&  (  2) &  1568768&  (  3) \\
  1570816&  (  2) &  1572864&  (384) &  1574912&  (  2) &  1576960&  (  4) &  1579008&  (  2) \\
  1581056&  (  7) &  1583104&  (  2) &  1585152&  (  4) &  1587200&  (  2) &  1589248&  ( 13) \\
  1591296&  (  2) &  1593344&  (  4) &  1595392&  (  2) &  1597440&  (  7) &  1599488&  (  2) \\
  1601536&  (  4) &  1603584&  (  2) &  1605632&  ( 25) &  1607680&  (  2) &  1609728&  (  4) \\
  1611776&  (  2) &  1613824&  (  7) &  1615872&  (  2) &  1617920&  (  4) &  1619968&  (  2) \\
  1622016&  ( 13) &  1624064&  (  2) &  1626112&  (  4) &  1628160&  (  2) &  1630208&  (  7) \\
  1632256&  (  2) &  1634304&  (  4) &  1636352&  (  2) &  1638400&  ( 50) &  1640448&  (  2) \\
  1642496&  (  4) &  1644544&  (  2) &  1646592&  (  7) &  1648640&  (  2) &  1650688&  (  4) \\
  1652736&  (  2) &  1654784&  ( 13) &  1656832&  (  2) &  1658880&  (  4) &  1660928&  (  2) \\
  1662976&  (  7) &  1665024&  (  2) &  1667072&  (  4) &  1669120&  (  2) &  1671168&  ( 26) \\
  1673216&  (  2) &  1675264&  (  4) &  1677312&  (  2) &  1679360&  (  7) &  1681408&  (  2) \\
  1683456&  (  4) &  1685504&  (  2) &  1687552&  ( 13) &  1689600&  (  2) &  1691648&  (  4) \\
  1693696&  (  2) &  1695744&  (  7) &  1697792&  (  2) &  1699840&  (  4) &  1701888&  (  2) \\
  1703936&  (104) &  1705984&  (  2) &  1708032&  (  4) &  1710080&  (  2) &  1712128&  (  7) \\
  1714176&  (  2) &  1716224&  (  4) &  1718272&  (  2) &  1720320&  ( 14) &  1722368&  (  2) \\
  1724416&  (  4) &  1726464&  (  2) &  1728512&  (  7) &  1730560&  (  2) &  1732608&  (  4) \\
  1734656&  (  2) &  1736704&  ( 27) &  1738752&  (  2) &  1740800&  (  4) &  1742848&  (  2) \\
  1744896&  (  7) &  1746944&  (  2) &  1748992&  (  4) &  1751040&  (  2) &  1753088&  ( 14) \\
  1755136&  (  2) &  1757184&  (  4) &  1759232&  (  2) &  1761280&  (  7) &  1763328&  (  2) \\
\end{tabular}
\newpage
\small
\begin{tabular}{|rr|rr|rr|rr|rr|}
  1765376&  (  4) &  1767424&  (  2) &  1769472&  ( 54) &  1771520&  (  2) &  1773568&  (  4) \\
  1775616&  (  2) &  1777664&  (  7) &  1779712&  (  2) &  1781760&  (  4) &  1783808&  (  2) \\
  1785856&  ( 14) &  1787904&  (  2) &  1789952&  (  4) &  1792000&  (  2) &  1794048&  (  7) \\
  1796096&  (  2) &  1798144&  (  4) &  1800192&  (  2) &  1802240&  ( 28) &  1804288&  (  2) \\
  1806336&  (  4) &  1808384&  (  2) &  1810432&  (  7) &  1812480&  (  2) &  1814528&  (  4) \\
  1816576&  (  2) &  1818624&  ( 14) &  1820672&  (  2) &  1822720&  (  4) &  1824768&  (  2) \\
  1826816&  (  7) &  1828864&  (  2) &  1830912&  (  4) &  1832960&  (  2) &  1835008&  (224) \\
  1837056&  (  2) &  1839104&  (  4) &  1841152&  (  2) &  1843200&  (  8) &  1845248&  (  2) \\
  1847296&  (  4) &  1849344&  (  2) &  1851392&  ( 15) &  1853440&  (  2) &  1855488&  (  4) \\
  1857536&  (  2) &  1859584&  (  8) &  1861632&  (  2) &  1863680&  (  4) &  1865728&  (  2) \\
  1867776&  ( 29) &  1869824&  (  2) &  1871872&  (  4) &  1873920&  (  2) &  1875968&  (  8) \\
  1878016&  (  2) &  1880064&  (  4) &  1882112&  (  2) &  1884160&  ( 15) &  1886208&  (  2) \\
  1888256&  (  4) &  1890304&  (  2) &  1892352&  (  8) &  1894400&  (  2) &  1896448&  (  4) \\
  1898496&  (  2) &  1900544&  ( 58) &  1902592&  (  2) &  1904640&  (  4) &  1906688&  (  2) \\
  1908736&  (  8) &  1910784&  (  2) &  1912832&  (  4) &  1914880&  (  2) &  1916928&  ( 15) \\
  1918976&  (  2) &  1921024&  (  4) &  1923072&  (  2) &  1925120&  (  8) &  1927168&  (  2) \\
  1929216&  (  4) &  1931264&  (  2) &  1933312&  ( 30) &  1935360&  (  2) &  1937408&  (  4) \\
  1939456&  (  2) &  1941504&  (  8) &  1943552&  (  2) &  1945600&  (  4) &  1947648&  (  2) \\
  1949696&  ( 15) &  1951744&  (  2) &  1953792&  (  4) &  1955840&  (  2) &  1957888&  (  8) \\
  1959936&  (  2) &  1961984&  (  4) &  1964032&  (  2) &  1966080&  (120) &  1968128&  (  2) \\
  1970176&  (  4) &  1972224&  (  2) &  1974272&  (  8) &  1976320&  (  2) &  1978368&  (  4) \\
  1980416&  (  2) &  1982464&  ( 16) &  1984512&  (  2) &  1986560&  (  4) &  1988608&  (  2) \\
  1990656&  (  8) &  1992704&  (  2) &  1994752&  (  4) &  1996800&  (  2) &  1998848&  ( 31) \\
  2000896&  (  2) &  2002944&  (  4) &  2004992&  (  2) &  2007040&  (  8) &  2009088&  (  2) \\
  2011136&  (  4) &  2013184&  (  2) &  2015232&  ( 16) &  2017280&  (  2) &  2019328&  (  4) \\
  2021376&  (  2) &  2023424&  (  8) &  2025472&  (  2) &  2027520&  (  4) &  2029568&  (  2) \\
  2031616&  ( 62) &  2033664&  (  2) &  2035712&  (  4) &  2037760&  (  2) &  2039808&  (  8) \\
  2041856&  (  2) &  2043904&  (  4) &  2045952&  (  2) &  2048000&  ( 16) &  2050048&  (  2) \\
  2052096&  (  4) &  2054144&  (  2) &  2056192&  (  8) &  2058240&  (  2) &  2060288&  (  4) \\
  2062336&  (  2) &  2064384&  ( 32) &  2066432&  (  2) &  2068480&  (  4) &  2070528&  (  2) \\
  2072576&  (  8) &  2074624&  (  2) &  2076672&  (  4) &  2078720&  (  2) &  2080768&  ( 16) \\
  2082816&  (  2) &  2084864&  (  4) &  2086912&  (  2) &  2088960&  (  8) &  2091008&  (  2) \\
  2093056&  (  4) &  2095104&  (  2) &  2097152& (2048) &  2098176&  (  2) &  2099200&  (  3) \\
  2100224&  (  2) &  2101248&  (  5) &  2102272&  (  2) &  2103296&  (  3) &  2104320&  (  2) \\
  2105344&  (  9) &  2106368&  (  2) &  2107392&  (  3) &  2108416&  (  2) &  2109440&  (  5) \\
  2110464&  (  2) &  2111488&  (  3) &  2112512&  (  2) &  2113536&  ( 17) &  2114560&  (  2) \\
  2115584&  (  3) &  2116608&  (  2) &  2117632&  (  5) &  2118656&  (  2) &  2119680&  (  3) \\
  2120704&  (  2) &  2121728&  (  9) &  2122752&  (  2) &  2123776&  (  3) &  2124800&  (  2) \\
  2125824&  (  5) &  2126848&  (  2) &  2127872&  (  3) &  2128896&  (  2) &  2129920&  ( 33) \\
  2130944&  (  2) &  2131968&  (  3) &  2132992&  (  2) &  2134016&  (  5) &  2135040&  (  2) \\
  2136064&  (  3) &  2137088&  (  2) &  2138112&  (  9) &  2139136&  (  2) &  2140160&  (  3) \\
  2141184&  (  2) &  2142208&  (  5) &  2143232&  (  2) &  2144256&  (  3) &  2145280&  (  2) \\
  2146304&  ( 17) &  2147328&  (  2) &  2148352&  (  3) &  2149376&  (  2) &  2150400&  (  5) \\
  2151424&  (  2) &  2152448&  (  3) &  2153472&  (  2) &  2154496&  (  9) &  2155520&  (  2) \\
  2156544&  (  3) &  2157568&  (  2) &  2158592&  (  5) &  2159616&  (  2) &  2160640&  (  3) \\
  2161664&  (  2) &  2162688&  ( 66) &  2163712&  (  2) &  2164736&  (  3) &  2165760&  (  2) \\
  2166784&  (  5) &  2167808&  (  2) &  2168832&  (  3) &  2169856&  (  2) &  2170880&  (  9) \\
  2171904&  (  2) &  2172928&  (  3) &  2173952&  (  2) &  2174976&  (  5) &  2176000&  (  2) \\
  2177024&  (  3) &  2178048&  (  2) &  2179072&  ( 17) &  2180096&  (  2) &  2181120&  (  3) \\
  2182144&  (  2) &  2183168&  (  5) &  2184192&  (  2) &  2185216&  (  3) &  2186240&  (  2) \\
  2187264&  (  9) &  2188288&  (  2) &  2189312&  (  3) &  2190336&  (  2) &  2191360&  (  5) \\
  2192384&  (  2) &  2193408&  (  3) &  2194432&  (  2) &  2195456&  ( 34) &  2196480&  (  2) \\
  2197504&  (  3) &  2198528&  (  2) &  2199552&  (  5) &  2200576&  (  2) &  2201600&  (  3) \\
\end{tabular}
\newpage
\small
\begin{tabular}{|rr|rr|rr|rr|rr|}
  2202624&  (  2) &  2203648&  (  9) &  2204672&  (  2) &  2205696&  (  3) &  2206720&  (  2) \\
  2207744&  (  5) &  2208768&  (  2) &  2209792&  (  3) &  2210816&  (  2) &  2211840&  ( 17) \\
  2212864&  (  2) &  2213888&  (  3) &  2214912&  (  2) &  2215936&  (  5) &  2216960&  (  2) \\
  2217984&  (  3) &  2219008&  (  2) &  2220032&  (  9) &  2221056&  (  2) &  2222080&  (  3) \\
  2223104&  (  2) &  2224128&  (  5) &  2225152&  (  2) &  2226176&  (  3) &  2227200&  (  2) \\
  2228224&  (136) &  2229248&  (  2) &  2230272&  (  3) &  2231296&  (  2) &  2232320&  (  5) \\
  2233344&  (  2) &  2234368&  (  3) &  2235392&  (  2) &  2236416&  (  9) &  2237440&  (  2) \\
  2238464&  (  3) &  2239488&  (  2) &  2240512&  (  5) &  2241536&  (  2) &  2242560&  (  3) \\
  2243584&  (  2) &  2244608&  ( 18) &  2245632&  (  2) &  2246656&  (  3) &  2247680&  (  2) \\
  2248704&  (  5) &  2249728&  (  2) &  2250752&  (  3) &  2251776&  (  2) &  2252800&  (  9) \\
  2253824&  (  2) &  2254848&  (  3) &  2255872&  (  2) &  2256896&  (  5) &  2257920&  (  2) \\
  2258944&  (  3) &  2259968&  (  2) &  2260992&  ( 35) &  2262016&  (  2) &  2263040&  (  3) \\
  2264064&  (  2) &  2265088&  (  5) &  2266112&  (  2) &  2267136&  (  3) &  2268160&  (  2) \\
  2269184&  (  9) &  2270208&  (  2) &  2271232&  (  3) &  2272256&  (  2) &  2273280&  (  5) \\
  2274304&  (  2) &  2275328&  (  3) &  2276352&  (  2) &  2277376&  ( 18) &  2278400&  (  2) \\
  2279424&  (  3) &  2280448&  (  2) &  2281472&  (  5) &  2282496&  (  2) &  2283520&  (  3) \\
  2284544&  (  2) &  2285568&  (  9) &  2286592&  (  2) &  2287616&  (  3) &  2288640&  (  2) \\
  2289664&  (  5) &  2290688&  (  2) &  2291712&  (  3) &  2292736&  (  2) &  2293760&  ( 70) \\
  2294784&  (  2) &  2295808&  (  3) &  2296832&  (  2) &  2297856&  (  5) &  2298880&  (  2) \\
  2299904&  (  3) &  2300928&  (  2) &  2301952&  (  9) &  2302976&  (  2) &  2304000&  (  3) \\
  2305024&  (  2) &  2306048&  (  5) &  2307072&  (  2) &  2308096&  (  3) &  2309120&  (  2) \\
  2310144&  ( 18) &  2311168&  (  2) &  2312192&  (  3) &  2313216&  (  2) &  2314240&  (  5) \\
  2315264&  (  2) &  2316288&  (  3) &  2317312&  (  2) &  2318336&  (  9) &  2319360&  (  2) \\
  2320384&  (  3) &  2321408&  (  2) &  2322432&  (  5) &  2323456&  (  2) &  2324480&  (  3) \\
  2325504&  (  2) &  2326528&  ( 36) &  2327552&  (  2) &  2328576&  (  3) &  2329600&  (  2) \\
  2330624&  (  5) &  2331648&  (  2) &  2332672&  (  3) &  2333696&  (  2) &  2334720&  (  9) \\
  2335744&  (  2) &  2336768&  (  3) &  2337792&  (  2) &  2338816&  (  5) &  2339840&  (  2) \\
  2340864&  (  3) &  2341888&  (  2) &  2342912&  ( 18) &  2343936&  (  2) &  2344960&  (  3) \\
  2345984&  (  2) &  2347008&  (  5) &  2348032&  (  2) &  2349056&  (  3) &  2350080&  (  2) \\
  2351104&  (  9) &  2352128&  (  2) &  2353152&  (  3) &  2354176&  (  2) &  2355200&  (  5) \\
  2356224&  (  2) &  2357248&  (  3) &  2358272&  (  2) &  2359296&  (288) &  2360320&  (  2) \\
  2361344&  (  3) &  2362368&  (  2) &  2363392&  (  5) &  2364416&  (  2) &  2365440&  (  3) \\
  2366464&  (  2) &  2367488&  ( 10) &  2368512&  (  2) &  2369536&  (  3) &  2370560&  (  2) \\
  2371584&  (  5) &  2372608&  (  2) &  2373632&  (  3) &  2374656&  (  2) &  2375680&  ( 19) \\
  2376704&  (  2) &  2377728&  (  3) &  2378752&  (  2) &  2379776&  (  5) &  2380800&  (  2) \\
  2381824&  (  3) &  2382848&  (  2) &  2383872&  ( 10) &  2384896&  (  2) &  2385920&  (  3) \\
  2386944&  (  2) &  2387968&  (  5) &  2388992&  (  2) &  2390016&  (  3) &  2391040&  (  2) \\
  2392064&  ( 37) &  2393088&  (  2) &  2394112&  (  3) &  2395136&  (  2) &  2396160&  (  5) \\
  2397184&  (  2) &  2398208&  (  3) &  2399232&  (  2) &  2400256&  ( 10) &  2401280&  (  2) \\
  2402304&  (  3) &  2403328&  (  2) &  2404352&  (  5) &  2405376&  (  2) &  2406400&  (  3) \\
  2407424&  (  2) &  2408448&  ( 19) &  2409472&  (  2) &  2410496&  (  3) &  2411520&  (  2) \\
  2412544&  (  5) &  2413568&  (  2) &  2414592&  (  3) &  2415616&  (  2) &  2416640&  ( 10) \\
  2417664&  (  2) &  2418688&  (  3) &  2419712&  (  2) &  2420736&  (  5) &  2421760&  (  2) \\
  2422784&  (  3) &  2423808&  (  2) &  2424832&  ( 74) &  2425856&  (  2) &  2426880&  (  3) \\
  2427904&  (  2) &  2428928&  (  5) &  2429952&  (  2) &  2430976&  (  3) &  2432000&  (  2) \\
  2433024&  ( 10) &  2434048&  (  2) &  2435072&  (  3) &  2436096&  (  2) &  2437120&  (  5) \\
  2438144&  (  2) &  2439168&  (  3) &  2440192&  (  2) &  2441216&  ( 19) &  2442240&  (  2) \\
  2443264&  (  3) &  2444288&  (  2) &  2445312&  (  5) &  2446336&  (  2) &  2447360&  (  3) \\
  2448384&  (  2) &  2449408&  ( 10) &  2450432&  (  2) &  2451456&  (  3) &  2452480&  (  2) \\
  2453504&  (  5) &  2454528&  (  2) &  2455552&  (  3) &  2456576&  (  2) &  2457600&  ( 38) \\
  2458624&  (  2) &  2459648&  (  3) &  2460672&  (  2) &  2461696&  (  5) &  2462720&  (  2) \\
  2463744&  (  3) &  2464768&  (  2) &  2465792&  ( 10) &  2466816&  (  2) &  2467840&  (  3) \\
  2468864&  (  2) &  2469888&  (  5) &  2470912&  (  2) &  2471936&  (  3) &  2472960&  (  2) \\
\end{tabular}
\newpage
\small
\begin{tabular}{|rr|rr|rr|rr|rr|}
  2473984&  ( 19) &  2475008&  (  2) &  2476032&  (  3) &  2477056&  (  2) &  2478080&  (  5) \\
  2479104&  (  2) &  2480128&  (  3) &  2481152&  (  2) &  2482176&  ( 10) &  2483200&  (  2) \\
  2484224&  (  3) &  2485248&  (  2) &  2486272&  (  5) &  2487296&  (  2) &  2488320&  (  3) \\
  2489344&  (  2) &  2490368&  (152) &  2491392&  (  2) &  2492416&  (  3) &  2493440&  (  2) \\
  2494464&  (  5) &  2495488&  (  2) &  2496512&  (  3) &  2497536&  (  2) &  2498560&  ( 10) \\
  2499584&  (  2) &  2500608&  (  3) &  2501632&  (  2) &  2502656&  (  5) &  2503680&  (  2) \\
  2504704&  (  3) &  2505728&  (  2) &  2506752&  ( 20) &  2507776&  (  2) &  2508800&  (  3) \\
  2509824&  (  2) &  2510848&  (  5) &  2511872&  (  2) &  2512896&  (  3) &  2513920&  (  2) \\
  2514944&  ( 10) &  2515968&  (  2) &  2516992&  (  3) &  2518016&  (  2) &  2519040&  (  5) \\
  2520064&  (  2) &  2521088&  (  3) &  2522112&  (  2) &  2523136&  ( 39) &  2524160&  (  2) \\
  2525184&  (  3) &  2526208&  (  2) &  2527232&  (  5) &  2528256&  (  2) &  2529280&  (  3) \\
  2530304&  (  2) &  2531328&  ( 10) &  2532352&  (  2) &  2533376&  (  3) &  2534400&  (  2) \\
  2535424&  (  5) &  2536448&  (  2) &  2537472&  (  3) &  2538496&  (  2) &  2539520&  ( 20) \\
  2540544&  (  2) &  2541568&  (  3) &  2542592&  (  2) &  2543616&  (  5) &  2544640&  (  2) \\
  2545664&  (  3) &  2546688&  (  2) &  2547712&  ( 10) &  2548736&  (  2) &  2549760&  (  3) \\
  2550784&  (  2) &  2551808&  (  5) &  2552832&  (  2) &  2553856&  (  3) &  2554880&  (  2) \\
  2555904&  ( 78) &  2556928&  (  2) &  2557952&  (  3) &  2558976&  (  2) &  2560000&  (  5) \\
  2561024&  (  2) &  2562048&  (  3) &  2563072&  (  2) &  2564096&  ( 10) &  2565120&  (  2) \\
  2566144&  (  3) &  2567168&  (  2) &  2568192&  (  5) &  2569216&  (  2) &  2570240&  (  3) \\
  2571264&  (  2) &  2572288&  ( 20) &  2573312&  (  2) &  2574336&  (  3) &  2575360&  (  2) \\
  2576384&  (  5) &  2577408&  (  2) &  2578432&  (  3) &  2579456&  (  2) &  2580480&  ( 10) \\
  2581504&  (  2) &  2582528&  (  3) &  2583552&  (  2) &  2584576&  (  5) &  2585600&  (  2) \\
  2586624&  (  3) &  2587648&  (  2) &  2588672&  ( 40) &  2589696&  (  2) &  2590720&  (  3) \\
  2591744&  (  2) &  2592768&  (  5) &  2593792&  (  2) &  2594816&  (  3) &  2595840&  (  2) \\
  2596864&  ( 10) &  2597888&  (  2) &  2598912&  (  3) &  2599936&  (  2) &  2600960&  (  5) \\
  2601984&  (  2) &  2603008&  (  3) &  2604032&  (  2) &  2605056&  ( 20) &  2606080&  (  2) \\
  2607104&  (  3) &  2608128&  (  2) &  2609152&  (  5) &  2610176&  (  2) &  2611200&  (  3) \\
  2612224&  (  2) &  2613248&  ( 10) &  2614272&  (  2) &  2615296&  (  3) &  2616320&  (  2) \\
  2617344&  (  5) &  2618368&  (  2) &  2619392&  (  3) &  2620416&  (  2) &  2621440&  (640) \\
  2622464&  (  2) &  2623488&  (  3) &  2624512&  (  2) &  2625536&  (  6) &  2626560&  (  2) \\
  2627584&  (  3) &  2628608&  (  2) &  2629632&  ( 11) &  2630656&  (  2) &  2631680&  (  3) \\
  2632704&  (  2) &  2633728&  (  6) &  2634752&  (  2) &  2635776&  (  3) &  2636800&  (  2) \\
  2637824&  ( 21) &  2638848&  (  2) &  2639872&  (  3) &  2640896&  (  2) &  2641920&  (  6) \\
  2642944&  (  2) &  2643968&  (  3) &  2644992&  (  2) &  2646016&  ( 11) &  2647040&  (  2) \\
  2648064&  (  3) &  2649088&  (  2) &  2650112&  (  6) &  2651136&  (  2) &  2652160&  (  3) \\
  2653184&  (  2) &  2654208&  ( 41) &  2655232&  (  2) &  2656256&  (  3) &  2657280&  (  2) \\
  2658304&  (  6) &  2659328&  (  2) &  2660352&  (  3) &  2661376&  (  2) &  2662400&  ( 11) \\
  2663424&  (  2) &  2664448&  (  3) &  2665472&  (  2) &  2666496&  (  6) &  2667520&  (  2) \\
  2668544&  (  3) &  2669568&  (  2) &  2670592&  ( 21) &  2671616&  (  2) &  2672640&  (  3) \\
  2673664&  (  2) &  2674688&  (  6) &  2675712&  (  2) &  2676736&  (  3) &  2677760&  (  2) \\
  2678784&  ( 11) &  2679808&  (  2) &  2680832&  (  3) &  2681856&  (  2) &  2682880&  (  6) \\
  2683904&  (  2) &  2684928&  (  3) &  2685952&  (  2) &  2686976&  ( 82) &  2688000&  (  2) \\
  2689024&  (  3) &  2690048&  (  2) &  2691072&  (  6) &  2692096&  (  2) &  2693120&  (  3) \\
  2694144&  (  2) &  2695168&  ( 11) &  2696192&  (  2) &  2697216&  (  3) &  2698240&  (  2) \\
  2699264&  (  6) &  2700288&  (  2) &  2701312&  (  3) &  2702336&  (  2) &  2703360&  ( 21) \\
  2704384&  (  2) &  2705408&  (  3) &  2706432&  (  2) &  2707456&  (  6) &  2708480&  (  2) \\
  2709504&  (  3) &  2710528&  (  2) &  2711552&  ( 11) &  2712576&  (  2) &  2713600&  (  3) \\
  2714624&  (  2) &  2715648&  (  6) &  2716672&  (  2) &  2717696&  (  3) &  2718720&  (  2) \\
  2719744&  ( 42) &  2720768&  (  2) &  2721792&  (  3) &  2722816&  (  2) &  2723840&  (  6) \\
  2724864&  (  2) &  2725888&  (  3) &  2726912&  (  2) &  2727936&  ( 11) &  2728960&  (  2) \\
  2729984&  (  3) &  2731008&  (  2) &  2732032&  (  6) &  2733056&  (  2) &  2734080&  (  3) \\
  2735104&  (  2) &  2736128&  ( 21) &  2737152&  (  2) &  2738176&  (  3) &  2739200&  (  2) \\
  2740224&  (  6) &  2741248&  (  2) &  2742272&  (  3) &  2743296&  (  2) &  2744320&  ( 11) \\
\end{tabular}
\newpage
\small
\begin{tabular}{|rr|rr|rr|rr|rr|}
  2745344&  (  2) &  2746368&  (  3) &  2747392&  (  2) &  2748416&  (  6) &  2749440&  (  2) \\
  2750464&  (  3) &  2751488&  (  2) &  2752512&  (168) &  2753536&  (  2) &  2754560&  (  3) \\
  2755584&  (  2) &  2756608&  (  6) &  2757632&  (  2) &  2758656&  (  3) &  2759680&  (  2) \\
  2760704&  ( 11) &  2761728&  (  2) &  2762752&  (  3) &  2763776&  (  2) &  2764800&  (  6) \\
  2765824&  (  2) &  2766848&  (  3) &  2767872&  (  2) &  2768896&  ( 22) &  2769920&  (  2) \\
  2770944&  (  3) &  2771968&  (  2) &  2772992&  (  6) &  2774016&  (  2) &  2775040&  (  3) \\
  2776064&  (  2) &  2777088&  ( 11) &  2778112&  (  2) &  2779136&  (  3) &  2780160&  (  2) \\
  2781184&  (  6) &  2782208&  (  2) &  2783232&  (  3) &  2784256&  (  2) &  2785280&  ( 43) \\
  2786304&  (  2) &  2787328&  (  3) &  2788352&  (  2) &  2789376&  (  6) &  2790400&  (  2) \\
  2791424&  (  3) &  2792448&  (  2) &  2793472&  ( 11) &  2794496&  (  2) &  2795520&  (  3) \\
  2796544&  (  2) &  2797568&  (  6) &  2798592&  (  2) &  2799616&  (  3) &  2800640&  (  2) \\
  2801664&  ( 22) &  2802688&  (  2) &  2803712&  (  3) &  2804736&  (  2) &  2805760&  (  6) \\
  2806784&  (  2) &  2807808&  (  3) &  2808832&  (  2) &  2809856&  ( 11) &  2810880&  (  2) \\
  2811904&  (  3) &  2812928&  (  2) &  2813952&  (  6) &  2814976&  (  2) &  2816000&  (  3) \\
  2817024&  (  2) &  2818048&  ( 86) &  2819072&  (  2) &  2820096&  (  3) &  2821120&  (  2) \\
  2822144&  (  6) &  2823168&  (  2) &  2824192&  (  3) &  2825216&  (  2) &  2826240&  ( 11) \\
  2827264&  (  2) &  2828288&  (  3) &  2829312&  (  2) &  2830336&  (  6) &  2831360&  (  2) \\
  2832384&  (  3) &  2833408&  (  2) &  2834432&  ( 22) &  2835456&  (  2) &  2836480&  (  3) \\
  2837504&  (  2) &  2838528&  (  6) &  2839552&  (  2) &  2840576&  (  3) &  2841600&  (  2) \\
  2842624&  ( 11) &  2843648&  (  2) &  2844672&  (  3) &  2845696&  (  2) &  2846720&  (  6) \\
  2847744&  (  2) &  2848768&  (  3) &  2849792&  (  2) &  2850816&  ( 44) &  2851840&  (  2) \\
  2852864&  (  3) &  2853888&  (  2) &  2854912&  (  6) &  2855936&  (  2) &  2856960&  (  3) \\
  2857984&  (  2) &  2859008&  ( 11) &  2860032&  (  2) &  2861056&  (  3) &  2862080&  (  2) \\
  2863104&  (  6) &  2864128&  (  2) &  2865152&  (  3) &  2866176&  (  2) &  2867200&  ( 22) \\
  2868224&  (  2) &  2869248&  (  3) &  2870272&  (  2) &  2871296&  (  6) &  2872320&  (  2) \\
  2873344&  (  3) &  2874368&  (  2) &  2875392&  ( 11) &  2876416&  (  2) &  2877440&  (  3) \\
  2878464&  (  2) &  2879488&  (  6) &  2880512&  (  2) &  2881536&  (  3) &  2882560&  (  2) \\
  2883584&  (352) &  2884608&  (  2) &  2885632&  (  3) &  2886656&  (  2) &  2887680&  (  6) \\
  2888704&  (  2) &  2889728&  (  3) &  2890752&  (  2) &  2891776&  ( 12) &  2892800&  (  2) \\
  2893824&  (  3) &  2894848&  (  2) &  2895872&  (  6) &  2896896&  (  2) &  2897920&  (  3) \\
  2898944&  (  2) &  2899968&  ( 23) &  2900992&  (  2) &  2902016&  (  3) &  2903040&  (  2) \\
  2904064&  (  6) &  2905088&  (  2) &  2906112&  (  3) &  2907136&  (  2) &  2908160&  ( 12) \\
  2909184&  (  2) &  2910208&  (  3) &  2911232&  (  2) &  2912256&  (  6) &  2913280&  (  2) \\
  2914304&  (  3) &  2915328&  (  2) &  2916352&  ( 45) &  2917376&  (  2) &  2918400&  (  3) \\
  2919424&  (  2) &  2920448&  (  6) &  2921472&  (  2) &  2922496&  (  3) &  2923520&  (  2) \\
  2924544&  ( 12) &  2925568&  (  2) &  2926592&  (  3) &  2927616&  (  2) &  2928640&  (  6) \\
  2929664&  (  2) &  2930688&  (  3) &  2931712&  (  2) &  2932736&  ( 23) &  2933760&  (  2) \\
  2934784&  (  3) &  2935808&  (  2) &  2936832&  (  6) &  2937856&  (  2) &  2938880&  (  3) \\
  2939904&  (  2) &  2940928&  ( 12) &  2941952&  (  2) &  2942976&  (  3) &  2944000&  (  2) \\
  2945024&  (  6) &  2946048&  (  2) &  2947072&  (  3) &  2948096&  (  2) &  2949120&  ( 90) \\
  2950144&  (  2) &  2951168&  (  3) &  2952192&  (  2) &  2953216&  (  6) &  2954240&  (  2) \\
  2955264&  (  3) &  2956288&  (  2) &  2957312&  ( 12) &  2958336&  (  2) &  2959360&  (  3) \\
  2960384&  (  2) &  2961408&  (  6) &  2962432&  (  2) &  2963456&  (  3) &  2964480&  (  2) \\
  2965504&  ( 23) &  2966528&  (  2) &  2967552&  (  3) &  2968576&  (  2) &  2969600&  (  6) \\
  2970624&  (  2) &  2971648&  (  3) &  2972672&  (  2) &  2973696&  ( 12) &  2974720&  (  2) \\
  2975744&  (  3) &  2976768&  (  2) &  2977792&  (  6) &  2978816&  (  2) &  2979840&  (  3) \\
  2980864&  (  2) &  2981888&  ( 46) &  2982912&  (  2) &  2983936&  (  3) &  2984960&  (  2) \\
  2985984&  (  6) &  2987008&  (  2) &  2988032&  (  3) &  2989056&  (  2) &  2990080&  ( 12) \\
  2991104&  (  2) &  2992128&  (  3) &  2993152&  (  2) &  2994176&  (  6) &  2995200&  (  2) \\
  2996224&  (  3) &  2997248&  (  2) &  2998272&  ( 23) &  2999296&  (  2) &         &
\end{tabular}

\end{document}